\documentclass[12pt,a4paper]{amsart}
\usepackage{amsfonts}
\usepackage{epsfig}
\usepackage{graphicx}
\usepackage{amsmath}
\usepackage{amssymb}
\usepackage{color}

\newcommand{\Rr}{\mathbb{R}}

\DeclareMathAlphabet{\mathpzc}{OT1}{pzc}{m}{it}

\newcounter{main}

\newtheorem{theorem}{Theorem}[section]

\newtheorem{lemma}[theorem]{Lemma}

\newtheorem{maintheorem}{Theorem}
\newtheorem{maincorollary}{Corollary}
\newtheorem{mainlemma}{Main Lemma}

\newtheorem{conjecture}{Conjecture}
\newtheorem{question}{Question}

\newcommand{\blanksquare}{\,\,\,$\sqcup\!\!\!\!\sqcap$}

\newcounter{example}
{{\stepcounter{example}}{\flushleft {\bf Example \arabic{example}:}}}%
{\par}

\title[Stably weakly shadowing and partial hyperbolicity]{Stably weakly shadowing symplectic maps are partially hyperbolic}

\author[M. Bessa]{M\'{a}rio Bessa}

\address{Universidade da Beira Interior, Rua Marqu\^es d'\'Avila e Bolama,
  6201-001 Covilh\~a
Portugal.}
\email{bessa@fc.up.pt}

\author[S. Vaz]{Sandra Vaz}

\address{Universidade da Beira Interior, Rua Marqu\^es d'\'Avila e Bolama,
  6201-001 Covilh\~a
Portugal.}
\email{svaz@ubi.pt}

\date{\today}

\begin{document}

\begin{abstract}
Let $M$ be a closed, symplectic connected Riemannian manifold, $f$ a symplectomorphism  on $M$. We prove that if $f$ is $C^1$-stably weakly shadowing on $M$, then the whole manifold $M$ admits a partially hyperbolic splitting. 
\end{abstract}

\maketitle

\section{Introduction, basic definitions and statement of the results}

\begin{subsection}{Introduction}
It is an old problem in smooth dynamics to perceive how the stability of a certain property on the phase space implies some hyperbolic-type of behavior of the tangent map of the system. The celebrated structural stability conjecture is the most important example of that (see e.g. ~\cite{Pu}). The $C^1$-stability of certain properties like topological conjugacy, topological stability, shadowing, expansiveness, specification, hyperbolic periodic orbits, etc, has been a debated issue in recent years. Here, we are interested in a shadowing-like property. 

The notion of shadowing applied to dynamics, introduced by Bowen \cite{B} in the context of hyperbolic dynamics, is motivated by the numerical computational idea of estimating differences between true and approximate solutions along orbits and to understand the influence of the errors that we commit and allow on each iterate.  In rough terms, we may ask if it is possible to obtain shadowing of approximate trajectories in a given dynamical system by true orbits of the system. We refer the reader to the monograph by Pilyugin's ~\cite{P} for a fairly complete exposition on the subject.

Several weakened shadowing-type definitions were developed in the last years, e.g. average shadowing \cite{S3}, asymptotic average shadowing \cite{G}, weak shadowing \cite{CP}. The weakly shadowing property first appear in the just mentioned paper by Corless and Pilyugin when related to the $C^0$-genericity of shadowing among dynamical systems. Informally speaking weakly shadowing allows that the pseudo-orbits may be approximated by true orbits if one forgets the time parametrization and consider only the distance between the orbit and the pseudo-orbit as two sets in the phase space. For examples of systems without the weak shadowing property see ~\cite[Example 2.12]{P} and for examples of systems satisfying the weakly shadowing property and not the shadowing one see ~\cite[Example 2.13]{P}. We intend in this paper to study the weakly shadowing property for diffeomorphisms which preserves a symplectic form. These systems are called \emph{symplectomorphisms}.

There are, of course, some limitations to the amount of information we can obtain from a given specific system that exhibits some shadowing-type property, since a $C^1$-close system may be absent of that property. For this reason it is of great utility and natural to consider that a selected model can be slightly perturbed in order to obtain the same property - the stably weakly shadowable dynamical systems. Nevertheless, it is worth to mention that $C^1$-stability in the symplectic setting only allows us to consider $C^1$-perturbations which are symplectic and not living in the broader space of diffeomorphisms or even in volume-preserving diffeomorphisms. This observation imply that the following results already proved for dissipative or volume-preserving diffeomorphisms are not applicable to our conservative context.

In ~\cite{S2} it is proved that if a diffeomorphism defined in a surface has the $C^1$-stable weak shadowing property, then it satisfy the axiom A and the no-cycle condition. However, the converse does not hold (see \cite{P}). We refer the paper \cite{PST} for more details on the relation between $C^1$-stability of weakly shadowing systems and structural stability in surfaces. Crovisier, in ~\cite{C}, proved that the weak shadowing property is generic (in the $C^1$-sense) for diffeomorphisms in closed manifolds. More precisely, he proved that for $C^1$-generic diffeomorphisms any pseudo-orbit is approximated in the Hausdorff topology by a finite segment of a genuine orbit. We observe that Crovisier results hold also in the conservative context (see ~\cite[\S 2.5]{C}).  Recently, Lee proved (see~\cite{L}) that $C^1$-weakly shadowing in area-preserving surfaces are Anosov. In this paper we generalize the results in \cite{GSW} and \cite{Y}. As a final remark we observe that our proof is strongly rooted on a result by Saghin and Xia \cite{SX} and holds only for the symplectic context. A different approach must be used for volume-preserving diffeomorphisms. However, in higher dimensions only a dominated splitting structure should be attainable. This is the subject of an ongoing work. 

\end{subsection}

\begin{subsection}{Basic definitions}

Denote by $M$ a $2d$-dimensional manifold ($d\geq 1$) with Riemaniann structure and endowed with a symplectic form  $\omega$. Let $\text{Diff}_\omega^{\,\,1}(M)$ denote the set of symplectomorphisms, i.e., of diffeomorphisms $f$ defined on $M$ and such that 
$$\omega_x(v_1,v_2)=\omega_{f(x)}(Df(x)\cdot v_1,Df(x)\cdot v_2),$$
for $x\in M$ and $v_1,v_2\in T_x M$. Consider this space endowed with the $C^1$ Whitney topology. It is well-known that $\text{Diff}_\omega^{\,\,1}(M)$ is a subset of all $C^1$ volume-preserving diffeomorphisms. The Riemannian inner-product induces a norm $\|\cdot\|$ on the tangent bundle $T_x M$. Denote the Riemannian distance by $d(\cdot,\cdot)$ and the open ball $B(x,r):=\{y\in M\colon d(x,y)<r\}$. We will use the usual uniform norm of a bounded linear map $A$ given by $\|A\|=\sup_{\|v\|=1}\|A\cdot v\|$.  By the theorem of Darboux (see e.g.~\cite[Theorem 1.18]{MZ}) there exists an atlas $\{\varphi_j^{-1}\colon U_j\to\Rr^{2n}\}$, where $U_j$ is an open subset of $M$, satisfying $\varphi_j^*\omega_0=\omega$ with $\omega_0=\sum_{i=1}^d dy_i\land dy_{d+i}$.

Fix some diffeomorphism $f\in \text{Diff}_\omega^{\,\,1}(M)$. Given $\delta>0$, we say that a sequence of points $\{x_i\}_{i\in\mathbb{Z}}\subset M$ is a \emph{$\delta$-pseudo-orbit of $f$} if $d(f(x_i),x_i)<\delta$ for all $i\in\mathbb{Z}$. We say that a sequence of points $\{x_i\}_{i\in\mathbb{Z}}\subset M$ is \emph{weakly $\epsilon$-shadowed by the $f$-orbit of $x$} if $\{x_i\}_{i\in\mathbb{Z}}\subset B({f^i(x),\epsilon})$ for all $i\in\mathbb{Z}$. For an $f$-invariant closed set $\Lambda\subseteq M$ we say that $f|_\Lambda$ has the \emph{weak shadowing property} if for every $\epsilon>0$, there exists $\delta>0$ such that for any $\delta$-pseudo-orbit $\{x_i\}_{i\in\mathbb{Z}}\subset \Lambda$ of $f$, there exists $x\in M$ such that the $f$-orbit of $x$ weakly $\epsilon$-shadows $\{x_i\}_{i\in\mathbb{Z}}$, i.e., $\{x_i\}_{i\in\mathbb{Z}}\subset B({f^i(x),\epsilon})$ for all $i\in\mathbb{Z}$. The diffeomorphism $f$ has the \emph{weak shadowing property} if $M=\Lambda$ in the above definition. We notice (cf. \cite{GSW}) that $f$ has the weak shadowing property if and only if $f^n$ has the weak shadowing property for all $n\in\mathbb{Z}$.

We say that $f$ is \emph{$C^1$-stably weakly shadowing} (in $M$) if there is a $C^1$-neighborhood $U(f)\cap \text{Diff}_\omega^{\,\,1}(M)$ of $f$, such that any symplectomorphism $g\in U(f)$ has the weakly shadowing property. We denote by $\mathcal{WS}(M)$ the open subset of $C^1$-stably weakly shadowing symplectomorphisms in $M$.

Recall, that a set $\Lambda\subseteq M$ is \emph{transitive} if it has a forward dense orbit. A diffeomorphism is transitive if $M$ is a transitive set for $f$. We observe that a transitive diffeomorphism has the weakly shadowing property. Thus, $C^1$-stable transitivity, implies the $C^1$-stable weakly shadowing property. As we already said weakly shadowing is $C^1$-generic, so let $\mathcal{R}_1\subset\text{Diff}_\omega^{\,\,1}(M)$ be this residual subset. Moreover, by ~\cite{ABC} there exists a $C^1$-residual set  $\mathcal{R}_2\subset\text{Diff}_\omega^{\,\,1}(M)$ such that any $f\in\mathcal{R}_2$ is transitive (actually, topologically mixing by \cite{AbC}). Thus, $\mathcal{R}=\mathcal{R}_1\cap \mathcal{R}_2$ is a residual subset where any $f\in\mathcal{R}$ is transitive and also weakly shadowable. Since $\text{Diff}_\omega^{\,\,1}(M)$ is a Baire space $\mathcal{R}$ is also $C^1$-dense.

Given an $f$-invariant set $\Lambda\subseteq M$ we say that $\Lambda$ is \emph{uniformly hyperbolic}  if the tangent vector bundle over $\Lambda$ splits into two $Df$-invariant subbundles $T\Lambda=E^u\oplus E^s$ such that $\|Df|_{E^s}\|\leq 1/2$ and $\|Df^{-1}|_{E^u}\|\leq 1/2$. When $\Lambda=M$ we say that $f$ is \emph{Anosov}. Clearly, there are lots of Anosov diffeomorphisms which are not symplectic. We say that an $f$-invariant set $\Lambda\subseteq M$ admits an \emph{$\ell$-dominated splitting} if there exists a continuous decomposition of the tangent bundle $T\Lambda$ into $Df$-invariant subbundles $E$ and $F$ such that $$\|Df^{\ell}(x)|_{F}\|.\|(Df^{\ell}(x)|_{E})^{-1}\|\leq 1/2,$$ 
in this case we say $E \succ_\ell F$ (i.e. $E$ $\ell$-dominates $F$). 
Finally, we say that an $f$-invariant set $\Lambda\subseteq M$ is \emph{uniformly partially hyperbolic}, if we have a splitting $E^s\oplus E^c\oplus E^u$ of $T\Lambda$ such that $E^s$ is uniformly contracting, $E^u$ is uniformly expanding, $E^c \succ E^s$ and  $E^u \succ E^c$. When $M$ is partially hyperbolic for $f$ we say that $f$ is a \emph{partially hyperbolic diffeomorphism}. It is proved in \cite{BV} that, in the symplectic world, the existence of a dominated splitting implies partial hyperbolicity. Let $\mathcal{PH}_\omega(M)\subset \text{Diff}_\omega^{\,\,1}(M)$ denote the $C^1$-open subset of partially hyperbolic symplectomorphisms.

A \emph{periodic orbit} for a diffeomorphism $f$ is a point $p\in M$ such that $f^\pi(p)=p$ where $\pi$ is the least positive integer satisfying the equality.  Given a periodic orbit $p$ of period $\pi$ of a diffeomorphism $f$ we say that $p$ is:
\begin{itemize}
\item \emph{hyperbolic} if $Df^{\pi}(p)$ has no norm one eigenvalues;
\item \emph{elliptic} if $Df^{\pi}(p)$ has only non-real eigenvalues of norm one;
\item \emph{$k$-elliptic} if there are precisely $k$ pairs of non-real eigenvalues of norm one and
\item \emph{parabolic} (or \emph{almost-elliptic}) if $Df^{\pi}(p)$ has only eigenvalues of norm one and, at least, a pair equal to $\{-1,1\}$.
\end{itemize}

\end{subsection}

\begin{subsection}{Statement of the results and some applications}
As we already said, here we develop the generalized versions of the results in \cite{GSW}, \cite{Y} and \cite{HT}. Our main result is the following.

\begin{maintheorem}\label{T1}
If $f \in  \text{Diff}_\omega^{\,\,1}(M)$ is $C^1$-stably weakly shadowing, then $f\in\mathcal{PH}_\omega(M)$.
\end{maintheorem}

Notice that in \cite{GSW,Y} it is considered a ``local" statement, i.e., $C^1$-stably weakly shadowing on a subset $\Lambda\subseteq M$, hypothesis of transitivity on $\Lambda$ and also $\Lambda$ being a homoclinic class. Here, we deal with $C^1$-stably weakly shadowing on the whole manifold. Observe that the transitivity on $M$ is $C^1$-generic \cite{ABC} and moreover there is only a single homoclinic class \cite{ABC}. Nevertheless, despite the fact that we treat the global statement we do not have to consider the transitivity hypothesis. Observe also that, under the $C^1$-robustly transitivity assumption on $M$, we obtain the partial hyperbolicity structure immediately by the results of \cite{HT}. 

We note that Bonatti and Viana \cite[\S 6.2]{BonV} build an open subset of partially hyperbolic (but not Anosov) transitive diffeomorphisms on $3$-dimensional manifolds (see also \cite{T}). We believe that our Theorem ~\ref{T1} is optimal at least for dimension $6$. The reason for that is based on the fact that Bonatti and Viana's construction should be possible to made symplectic but with three degrees of freedom. 

\begin{conjecture}
There exists non empty open subsets of $C^1$-stably weakly shadowing of symplectomorphisms of dimension $\geq 6$ such that $M$ is not uniformly hyperbolic for $f$ on that open sets. 
\end{conjecture}

We point out that we do not know if in dimension $4$ the hyperbolicity can be obtained. So, we ask the following:

\begin{question}
Are there Ma\~n\'e-Bonatti-Viana's examples (cf. \cite[\S 7.1]{BDV}) of symplectomorphisms on dimension four which are $C^1$-stably weakly shadowing (or $C^1$-robust transitive) but not uniformly hyperbolic?
\end{question}

Observe that a positive answer to previous questions will imply that $1$-elliptic periodic points may coexist with $C^1$-stably weakly shadowing (thus partial hyperbolicity) when $\dim(M)\geq 4$. Actually, by a result of Newhouse \cite{N} if the symplectomorphism is not Anosov then $1$-elliptic points can be created by arbitrary small $C^1$-perturbations of the symplectomorphism.

Next, we present an easy consequence of the previous theorem applied to surfaces. Its proof relies in the simple fact that, in the two-dimensional conservative case, dominated splitting is tantamount to hyperbolicity (see also ~\cite{L}). 

\begin{maincorollary}\label{C1}
Let $f\in  \text{Diff}_\omega^{\,\,1}(M)$ where $M$ is a surface. If $f$ is $C^1$-stably weakly shadowing, then $f$ is Anosov.
\end{maincorollary}

Although it is somewhat misleading, since we use results in \cite{HT} and a powerful dichotomy \cite{SX} to prove Theorem ~\ref{T1}, our results can be used to deduce the Horita and Tahzibi theorem \cite{HT}. This theorem states that $C^1$-robust transitivity implies partial hyperbolicity. Since $C^1$-robust transitivity implies $C^1$-stably weakly shadowing on $M$ we obtain, using Theorem ~\ref{T1}, that $M$ has a partially hyperbolic splitting. 

\begin{maincorollary}\label{C2}(Horita and Tahzibi ~\cite{HT})
Any $C^1$-robustly transitive symplectomorphism is partially hyperbolic.
\end{maincorollary}

Finally, we obtain some $C^1$-stable weakly shadowing ergodic map near a $C^1$-stable weakly shadowing. Recall that $f$ is said to be \emph{ergodic} if the only $f$-invariant subsets have zero or full volume.

\begin{maincorollary}\label{C4}
If $f\in  \text{Diff}_\omega^{\,\,1}(M)$ is $C^1$-stably weakly shadowing, then $f$ is $C^1$-approximated by a $C^1$-robust transitive partially hyperbolic symplectomorphism on $M$ which is also ergodic.
\end{maincorollary}

The proof is a straightforward application of powerful results on symplectomorphisms; by Theorem ~\ref{T1}, $M$ admits a partially hyperbolic splitting. Moreover, there exists a $C^1$ open and dense subset of $\mathcal{PH}_{\omega}(M)$ such that $f$ is transitive, for all $f$ on that subset. This is a corollary of the main theorems on $C^1$-denseness of accessibility by Dolgopyat and Wilkinson's \cite{DW} combined with a result of Brin \cite{B}. Finally, just recall \cite[Theorem A]{ABW} which says that ergodicity is $C^1$-generic for partially hyperbolic symplectomorphisms.  

\end{subsection}

\section{Proof of Theorem ~\ref{T1}}

\begin{subsection}{A main lemma}
The following symplectic version of  Frank's Lemma will be very useful to prove our main lemma.

\begin{lemma}(Horita and Tahzibi \cite[Lemma 5.1]{HT})\label{Franks}
Let $f \in \text{Diff}_\omega^{\,\,1}(M)$ and $U(f)$ be given. Then there are $\delta_0>0$ and $U_0(f)$ such that for any $g \in U_0(f)$, a finite set $\{ x_1, x_2,...,x_l \}$, a neighborhood U of $\{ x_1, x_2,..., x_l \}$ and symplectic maps $L_i: T_{x_i}M \rightarrow T_{g(x_i)}M$ satisfying $||L_i - Dg(x_i)|| < \delta_0$ for all $1 \leq i \leq l$ there are $\epsilon_0>0$ and $\tilde{g} \in U(f)$ such that
\begin{enumerate}
\item [a)] $\tilde{g}(x)=g(x)$ if $x \in M \backslash U$
\item [b)] $\tilde{g} (x)= \varphi_{g(x_i)} \circ L_i \circ \varphi^{-1}_{x_i}(x)$ if $x \in B(x_i,\epsilon_0)$ 
\end{enumerate}

\end{lemma}

Assertion b) implies that $\tilde{g}(x)=g(x)$ if $x \in \{ x_1, x_2,....,x_l \}$ and $D\tilde{g}(x_i) =L_i$. 
 
The following result is more or less the symplectic version of ~\cite[Lemma 3.2]{GSW} (see also \cite[Proposition 3.5]{L}). 
 
\begin{mainlemma}\label{ML}
Let $f\in \mathcal{WS}(M)$ and $U_0(f)$ be given by Lemma ~\ref{Franks} with respect to $U(f)$. Then, for any $g \in U_0(f)$, g does not contains elliptic points.
\end{mainlemma}

\begin{proof}
Let  $\dim M=2n$ and let us suppose that there is a symplectomorphim $g \in U_0(f)$ that have a periodic elliptic point $p$. To avoid notational complexity, we will suppose that $g(p)=p$, for $p$ an elliptic point. Then $Dg(p)$ has $n$ pairs of non-real eigenvalues: $|z_i|=|\bar{z}_i|=1$ , $i=1,...,n$  with $T_pM=E^{L_1}_p \oplus... \oplus E^{L_n}_p$ and $\dim(E_p^{L_i})=2, i=1,...,n$.

By Lemma ~\ref{Franks} there are $\epsilon_0>0$ and $\tilde{g} \in U(f) \subset \mathcal{WS}(M)$ such that $\tilde{g}(p)= g(p)=p$ and $\tilde{g}(x)= \varphi_{g(p)} \circ Dg(p) \circ \varphi_p^{-1}(x)$ if $x \in B(p,\epsilon_0)$. 

The next computations will be yield in $E^{L_1}_p(\epsilon_0)$ (i.e. a ball of radius $\epsilon_0$, centered in $0$ and inside $E^{L_1}_p$) without loss of generality.

With a $C^{1}$-small modification of the map $Dg(p)$ we may suppose that there is $\ell_1>0$ (the minimum number) such that $D g^{\ell_1}(p)\cdot v=v$ for any $v \in E^{L_1}_p (\epsilon_0) \cap \varphi^{-1}_p (B(p,\epsilon_0))$ by Lemma ~\ref{Franks}. Note that the restriction $\tilde{g}|_{\varphi_p(E^{L_i}_p(\epsilon_0)) \cap B(p,\epsilon_0)}$ of the map is a rotation, $i=2,...,n$.

Take $v_0 \in E^{L_1}_p(\epsilon_0)$ such that $\|v_0\| = \frac{\epsilon_0}{4}$ and set 

$$ \mathcal{I}_p = \varphi_p( \{ t.v_0: 1 \leq t \leq (1 + \epsilon_0/4)\}) \cap B(p,\epsilon_0).$$

Then, $\mathcal{I}_p$ is an arc such that

\begin{itemize}
\item $\tilde{g}^i(\mathcal{I}_p) \cap \tilde{g}^j(\mathcal{I}_p)= \emptyset$ if $0 \leq i \neq j \leq \ell_1 -1$;
\item $\tilde{g}^{\ell_1} (\mathcal{I}_p)=\mathcal{I}_p$ and $\tilde{g}^{\ell_1}_{|\mathcal{I}_p}$ is the identity map.
\end{itemize}

Choose $0<\epsilon_1 < \frac{\epsilon_0}{4}$ small enough such that 
$$B(\tilde{g}^i (\mathcal{I}_p),\epsilon_1) \cap B(\tilde{g}^j(\mathcal{I}_p),\epsilon_1) = \emptyset, \text{for all } 1 \leq i \neq j \leq \ell_1 -1.$$

Put $\epsilon= \epsilon_1$ and let $0 < \delta < \epsilon$ be the number given by the weak shadowing property of $\tilde{g}$.

Now, we are going to construct a $\delta$-pseudo-orbit $\{x_k \}_{k \in \mathbb{Z}}$ of $\tilde{g}$ in $\bigcup_{i=0}^{\ell_1-1} \tilde{g}^i(\mathcal{I}_p)$ which cannot be weakly $\epsilon$-shadowed by any $\tilde{g}$-true orbit of a point in M.

We find a finite sequence $\{ v_k \}_{k=0}^T$ in $\{ t.v_0: 1 \leq t \leq (1 + \frac{\epsilon_0}{4} )\}$ for some $T>0$, such that $v_T=(1+ \frac{\epsilon_0}{4}).v_0$ and $|v_k - v_{k+1}| < \delta$ for $0 \leq k \leq T-1$. Here, the $v_k$ are chosen such that if $v_k=t_k.v_0$, then $t_k < t_{k+1}$ for $0 \leq k \leq (T-1)$. Finally, define a $\delta$-pseudo-orbit $\{x_k \}_{k \in \mathbb{Z}}$ of $\tilde{g}$ by:

\begin{itemize}
\item
$x_k=\tilde{g}^k(\varphi_p(v_0))$ for $k<0$;
\item
$x_k=x_{i \ell_1 + j}=\tilde{g}^j(\varphi_p(v_i))$ for $0 \leq i \leq T-1$, $0 \leq j \leq \ell_1 -1$;
\item
$x_k=\tilde{g}^{k-\ell_1 T}(\varphi_p(v_T))$ for $k \geq \ell_1 T$
\end{itemize}

Since $\tilde{g} \in \mathcal{WS}(M)$, there is $y \in M$ weakly $\epsilon$-shadowing $\{ x_k \}_{k \in \mathbb{Z}}.$

The local structure of $\tilde{g}$ in a neighborhood of $\mathcal{I}_p$ in M is the direct product of the identity map $\tilde{g}|_{\mathcal{I}_p}$ by  rotations  $\tilde{g}|_{\varphi_p(E^{L_i}_p(\epsilon_0)) \cap B(p,\epsilon_0)}, i=2,...,n$.
Without loss of generality, we may assume that $y \in B(x_0,\epsilon)$.

If $y \in \mathcal{I}_p$ then taking $\tilde{h}= \tilde{g}^{\ell_1}$ since $\tilde{h}^i(y)=y$ for $i \in \mathbb{Z}$, and so $d(\tilde{h}^i (y), x_T)> \epsilon$ by choice of $\epsilon$, which is a contradiction.

If $y \notin \mathcal{I}_p$ and $y$ has components in  $\varphi_p(E^{L_i}_p(\epsilon_0)) \cap B(x_0,\epsilon_0)$ for  $i=2,...,n$ then by a small $C^1$-perturbation of $\tilde{g}|_{\varphi_p(E^{L_i}_p(\epsilon_0)) \cap B(p,\epsilon_0)}$, using  Lemma ~\ref{Franks} we may assume that there are $\ell_i$ (the minimum number) such that $Dg^{\ell_i}(p)\cdot v=v$ for any $v \in E^{L_i}_p (\epsilon_0) \cap \varphi^{-1}_p (B(p,\epsilon_0)), \ i=2,...,n$. Take $\kappa=\text{LCM}\{\ell_i, i=2,...,n\}$ and $\tilde{h}=\tilde{g}^{\kappa}$ (where LCM stands for the lowest common multiple). Then $\tilde{h}^i(y) \in B(x_0,\epsilon_0)$ and $d(\tilde{h}^i(y), x_T)>\epsilon$, by choice of $\epsilon$, for all $i \in \mathbb{Z}$, which is also a contradiction. This proves the lemma.
\end{proof} 

\end{subsection}

\begin{subsection}{End of the proof of Theorem~\ref{T1}}\label{proofs}

In \cite{SX} it was proved the following result which will be crucial to obtain our results:

\begin{theorem} (\cite[Theorem 2]{SX})\label{T4} 
If $f\notin\mathcal{PH}_\omega(M)$, then for any open
subset $U\subset M$ there exists an arbitrarily small $C^1$-perturbation $g$ of $f$ such that $g$ has
an elliptic periodic orbit through $U$.
\end{theorem}

Since, as we already said, the set of symplectomorphisms equipped with the $C^1$-topology is a Baire space we get that the residual subset of the previous theorem is also $C^1$-dense.

\begin{proof} (of Theorem ~\ref{T1}) The proof is by contradiction. Let us assume that $f\in\mathcal{WS}(M)$ and $f\notin\mathcal{PH}_\omega(M)$. There exists a neighborhood $U(f)$ of $f$ such that any $g\in U(f)$ is in $\mathcal{WS}(M)$. By Theorem~\ref{T4} for any open subset $U\subset M$ there exists a $C^1$-perturbation $g\in U(f)$ such that $g$ has an elliptic periodic point in $U$ and this contradicts the Main Lemma~\ref{ML}.
\end{proof}

\end{subsection}

\section*{Acknowledgements}

MB was partially supported by the FCT-Funda\c{c}\~ao para a Ci\^encia e a Tecnologia, project PTDC/MAT/099493/2008. SV  was partially supported by the FCT-Funda\c{c}\~ao para a Ci\^encia e a Tecnologia, project PEst-OE/MAT/UI0212/2011. The authors would like to thank Shaobo Gan and Kazuhiro Sakai for providing us their papers.

\end{document}